\documentclass[12pt,english]{amsart}
\usepackage[T1]{fontenc}
\usepackage[latin9]{inputenc}
\usepackage{amssymb}

\makeatletter
\usepackage[all]{xy}
\usepackage{color}
\usepackage{mathrsfs}
\usepackage{amsfonts}
\usepackage{amsthm}
\newcommand{\pes}[2]{\langle #1,#2\rangle}

\newtheorem{teo}{Theorem}

\newtheorem{prop}{Proposition}

\newtheorem{lem}{Lemma}

\newcommand{\N}{\mathbb{N}}

\newcommand{\C}{\mathbb{C}}
\newcommand{\R}{\mathbb{R}}

\newcommand{\bi}{\begin{itemize}}
\newcommand{\ei}{\end{itemize}}
\newcommand{\bd}{\begin{description}}
\newcommand{\ed}{\end{description}}
\newcommand{\beq}{\begin{equation}}
\newcommand{\eeq}{\end{equation}}
\newcommand{\beqn}{\begin{eqnarray}}
\newcommand{\eeqn}{\end{eqnarray}}
\newcommand{\beqna}{\begin{eqnarray*}}
\newcommand{\eeqna}{\end{eqnarray*}}
\newcommand{\mnr}{\mathscr{M}_n(\mathbb{\R})}
\newcommand{\mnc}{\mathscr{M}_n(\mathbb{\C})}

\newcommand{\mmnc}{\mathscr{M}_{k,p}(\mathbb{\C})}

\makeatother

\usepackage{babel}

\begin{document}
\fontsize{12pt}{17pt}\selectfont

\title{Average-case Perturbations and Smooth Condition Numbers}
\author{
Diego Armentano
}
\maketitle
\begin{flushleft}
Centro de Matem\'atica\\
 Facultad de Ciencias\\
Universidad de la Rep\'ublica\\
Calle Igua $4225$\\
 11400 Montevideo, Uruguay\\
diego@cmat.edu.uy
\end{flushleft}

\vspace{5pt}

\begin{abstract}
In this paper we define a new condition number adapted to directionally uniform perturbations. The definitions and theorems can be applied to a large class of  problems.
We show the relation with the classical condition number, and study some interesting examples.
\end{abstract}

\vspace{10pt}

{\footnotesize{\textit{Keywords:} Condition Numbers, Random Matrices, Systems of Linear Equations, Finding Kernels, Eigenvalue and Eigenvector Problems, Solving Polynomial Systems of Equations.

\textit{AMS subject classifications:} 65F35; 65Y05
}}

\vspace{5pt}

\section{Introduction}

Let $X$ and $Y$ be  two real (or complex) Riemannian manifolds of real dimensions $m$ and $n$  ($m\geq n$) respectively associated to some computational problem,
where $X$ is the space of \textit{inputs} and $Y$ is the space of \textit{outputs}.
Let $V\subset X\times Y$ be the \textit{solution variety}, i.e. the subset of pairs $(x,y)$ such
that $y$ is an output corresponding to the input $x$. Let $\pi_1:V\to X$ and $\pi_2:V\to Y$ be
the canonical projections, and $\Sigma\subset V$ be the set of critical points of the projection $\pi_1$.

 In case $\dim V=\dim X$, for each $(x,y)\in V\setminus \Sigma$
 there is a differentiable function
locally defined from some neighborhoods $U_x$ and $U_y$ of $x\in X$ and $y\in Y$ respectively,
namely $$G:=\pi_2\circ\pi_1^{-1}|_{U_x}:U_x\to U_y.$$

Let us denote by $\pes\cdot\cdot_{x}$ and $\pes\cdot\cdot_{y}$ the
Riemannian (or Hermitian) inner product in the tangent spaces $T_xX$ and $T_yY$ at $x$ and $y$ respectively.
The derivative $DG(x):T_xX\to T_yY$ is called the \textit{condition matrix} at $(x,y)$. The \textit{classical condition number} at $(x,y)\in V\setminus\Sigma$ is defined as
\beq\label{eq:CN}
\kappa(x,y):=
\max_{\substack{\dot{x}\in T_xX\\ {\|\dot{x}\|^2}_{x}=1}}\|DG(x)\dot{x}\|_{y}.
\eeq
This number is an upper-bound -to first-order approximation- of  the \textit{worst-case} sensitivity  of the output error with respect to small perturbations of the input. It plays an important role
to understand the behavior of algorithms and, as a consequence, appears in the usual bounds of the
running time of execution.
There exist an extensive literature about the role of the condition number in numerical analysis and complexity of algorithms, see for example \cite{Tref} and references therein.

In many practical situations,
there exist a discrepancy between theoretical analysis and observed performance of an algorithm.
There exist several approaches that attempt to rectify this discrepancy. Among them we find
\textit{average-case analysis} (see \cite{Edelman,Smale}) and \textit{smooth analysis} (see \cite{SpTeng,BCL,Wschebor}). For a comprehensive review on this subject with historical notes see \cite{Bur}.

In this paper, averaging is performed in a different form. In many problems, the space of inputs has a much
larger dimension than the one of the space of outputs ($m\gg n$). Then, it is natural to assume that infinitesimal
perturbations of the input will produce drastic changes in the output, only when they are performed in a few
directions. Then, a possibly different approach to analyze complexity of algorithms is to replace ``worst direction'' by a certain mean over all possible directions.
This alternative was already suggested and studied  in \cite{S-W-W-W} in the case of linear system solving $Ax=b$, and more generally, in \cite{Stewart} in the case of matrix perturbation theory where  the first-order perturbation expansion is assumed to be random.

In this paper we extend this approach to a large class of computational problems, restricting ourselves to the case of directionally uniform perturbations.

Generalizing the concept introduced in \cite{S-W-W-W} and \cite{Stewart}, we define the \textit{$p$th-average condition number} at $(x,y)$ as
\beq\label{eq:acn}
{\kappa_{av}}^{[p]}(x,y):=
\left[\frac{1}{\mbox{vol}(S_x^{m-1})}\int_{\dot{x}\in S_x^{m-1}}\|DG(x)\dot{x}\|_y^p\,dS_x^{m-1}(\dot x)  \right]^{1/p}
\quad (p=1,2,\ldots)
\eeq
where $\mbox{vol}(S_x^{m-1})=\frac{2\pi^{m/2}}{\Gamma(m/2)}$ is the
measure of the unit sphere $S_x^{m-1}$ in $T_xX$, and $dS_x^{m-1}$ is the induced volume element.
We will be mostly interested in the case $p=2$, which we simply write ${\kappa_{av}}$ and call it \textit{average condition number}.

Before the statement of the main theorem, we define the \textit{Frobenius condition number} as
 $$
\kappa_F(x,y):=\|DG(x)\|_F=\sqrt{\sigma_1^2+\cdots+\sigma_n^2}
$$
where $\|\cdot\|_F$ is the Frobenius norm and $\sigma_1,\ldots,\sigma_n$ are
the singular values of the condition matrix.

\begin{teo}\label{teo:average}
$$
{\kappa_{av}}^{[p]}(x,y)=
\frac{1}{\sqrt 2}\left[\frac{\Gamma\left(\frac{m}{2} \right)}{\Gamma\left(\frac{m+p}{2} \right)}\right]^{1/p}
\mathbb{E}(\|\eta_{\sigma_1,\ldots,\sigma_n}\|^p)^{1/p}
$$
where
$\|\cdot\|$ is the Euclidean norm in $\R^n$ and $\eta_{\sigma_1,\ldots,\sigma_n}$ is a centered Gaussian vector in $\R^n$ with diagonal covariance matrix
$Diag(\sigma_1^2,\ldots,\sigma_n^2)$.
\\
In particular,
\beq\label{teo:pdos}
{\kappa_{av}}(x,y)=\frac{\kappa_F(x,y)}{\sqrt m}.
\eeq
\end{teo}

This result is most interesting when $n\ll m$, for in that case
$$
{\kappa_{av}}(x,y)\leq \sqrt{\frac{n}{m}}\cdot\kappa(x,y)\ll \kappa(x,y),
$$
as $\kappa_F(x,y)\leq\sqrt{n}\cdot\kappa(x,y)$.
Thus, in these cases one may expect much better stability properties than those predicted by classical condition numbers.

In numerical analysis, many authors are interested in relative errors. Thus, in the case $(X,\|\cdot\|_X)$ and $(Y,\|\cdot\|_Y)$ are normed vector spaces, instead of consider the (absolute) condition number (\ref{eq:CN}), one could take the \textit{relative condition number} defined as
$$
\kappa_{rel}(x):= \frac{\|x \|_X}{\|y \|_Y}\kappa(x,y),
$$
and the \textit{relative Frobenius condition number} as
$$
{\kappa_{rel}}_F(x):=\frac{\|x \|_X}{\|y \|_Y} \kappa_F(x,y),
$$
(for simplicity we drop the $y$ in the argument).

In the same way, we define the \textit{relative $p$th-average condition number} as
\beq\label{def:cnrav}
{{\kappa }_{rel}}_{av}^{[p]}(x):=\frac{\|x \|_X}{\|y \|_Y}{ \kappa_{av}}^{[p]}(x,y),\quad (p=1,2,\ldots).
\eeq
For the case $p=2$ we simply write ${{\kappa }_{rel}}_{av}$ and call it \textit{relative average condition number}.

\emph{Theorem \ref{teo:average}} remains true if one change the (absolute) condition number by the relative condition number. In particular,
$$
{{\kappa }_{rel}}_{av}(x):=\frac{{{\kappa }_{rel}}_F(x,y)}{\sqrt m}.
$$

\vspace{5pt}

\section{Componentwise Analysis}

In the case $Y=\R^n$ we define the \textit{$k$th-componentwise condition number} at  $(x,y)\in V$ as:
\beq\label{eq:CNk}
\kappa(x,y,k):=\max_{\substack{\dot x\in T_xX\\{\|\dot x\|^2}_{x}=1}}|(DG(x)\dot x)_k|\qquad  (k=1,\ldots,n),
\eeq
where $|\cdot|$ is the absolute value and $w_k$ indicates the $k$th-component of the vector $w\in\R^n$.

Following \cite{S-W-W-W} for the linear case, we define the average componentwise condition number as
\beq\label{eq:cwacn}
{\kappa_{av}}^{[p]}(x,y,k):=
\left[\frac{1}{\mbox{vol}{(S_{x}^{m-1})}}\int_{\dot x\in S_{x}^{m-1}}\left|(DG(x)\dot x)_k\right|^p\,dS_{x}^{m-1}(\dot x)  \right]^{1/p}
\quad (p=1,2,\ldots).
\eeq
Then we have:
\begin{prop}\label{prop:averagecw}
$$
{\kappa_{av}}^{[p]}(x,y,k)=
\left[\frac{1}{\sqrt{\pi}}\cdot
\frac{\Gamma\left(\frac{m}{2} \right)}{\Gamma\left(\frac{m+p}{2} \right)}\cdot\Gamma\left(\frac{p+1}{2}\right)\right]^{1/p} \cdot{\kappa}_F(x,y,k).
$$
In particular,
$$
{\kappa_{av}}(x,y,k)=\frac{\kappa_F(x,y,k)}{\sqrt m} .
$$
\end{prop}

\vspace{10pt}

\begin{proof}
Observe that ${\kappa_{av}}^{[p]}(x,y,k)$ is the $p$th-average condition number for the problem
of finding the $k$th-component of $G=(G_1,\ldots,G_n)$.
\textit{Theorem \ref{teo:average}} applied to $G_k$ yields
$$
{\kappa_{av}}^{[p]}(x,y,k)=
\frac{1}{\sqrt 2}\left[\frac{\Gamma\left(\frac{m}{2} \right)}{\Gamma\left(\frac{m+p}{2} \right)}\right]^{\frac 1p}
\mathbb{E}(|\eta_{\sigma_1}|^p)^{1/p}
$$
where $\sigma_1=\|DG_k(x)\|=\kappa(x,y,k)$.
Then,
$$
\mathbb{E}(|\eta_{\sigma_1}|^p)^{1/p}=\kappa(x,y,k)\cdot \mathbb{E}(|\eta_1|^p)^{1/p}
$$
where $\eta_1$ is a standard normal in $\R$.
Finally,
$$
\mathbb{E}(|\eta_1|^p)=\frac{2}{\sqrt{2\pi}}\int_{0}^{\infty}\rho^p e^{-\rho^2/2}\,d\rho
=\frac{2}{\sqrt{2\pi}}2^{\frac{p-1}{2}}\Gamma(\frac{p+1}{2})
$$
and the proposition follows.
\end{proof}

\vspace{10pt}


\section{Examples}

In this section we will compute the average condition number for different problems:
systems of linear equations, eigenvalue and eigenvector problems, finding kernels of linear transformations and solving polynomial systems of equations.
The first two have been computed in \cite{Stewart} and are an easy consequence of \textit{Theorem \ref{teo:average}} and the usual condition number.
\vspace{10pt}

\subsection{Systems of Linear Equations}
We consider the problem of solving the system of linear equations $Ay=b$, where $A\in\mnr$ the space of $n\times n$ matrices with the Frobenius inner product, i.e. $\pes{A}{B}_F=\mbox{trace}(B^tA)$ ($B^t$ is the transpose of $B$), and $b\in\R^n$. \\
If we assume that $b$ is fixed, then the input space $X=\mnr$ with the Frobenius inner product, $Y=\R^n$ with the Euclidean inner product, and $\Sigma$ equals the subset of non-invertible matrices.
Then the map $G:\mnr\setminus\Sigma\to\R^n$ is globally defined and differentiable, namely
$$
G(A)=A^{-1}b \;(=y).
$$
By implicit differentiation,
\beq\label{eq:dglc}
DG(A)(\dot A)=-A^{-1}\dot A y.
\eeq
Is easy to see from (\ref{eq:dglc}) that $\kappa(A,y)=\|A^{-1}\|\cdot\|y\|$ and $\kappa_F(A,y)=\|A^{-1}\|_F\cdot\|y\|$. Then from \textit{Theorem \ref{teo:average}} we get
\beq\label{eq:acnlc}
{\kappa_{av}}(A,y)=\frac{\|A^{-1}\|_F\cdot\|y\|}{n}\leq
\frac{\kappa(A,y)}{\sqrt n}.
\eeq
A similar result was proved in \cite{Stewart}.

For the general case, we have $X=\mnr\times\R^n$ with the Frobenius inner product in $\mnr$ and the Euclidean inner product in $\R^n$. Then, $G:\mnr\setminus\Sigma\times\R^n\to\R^n$ satisfies $G(A,b)=A^{-1}b$.
\\
Is easy to see that $\kappa((A,b),y)=\|A^{-1}\|\cdot\sqrt{1+\|y\|^2}$ and $\kappa_F((A,b),y)=\|A^{-1}\|_F\cdot\sqrt{1+\|y\|^2}$.
Again from \textit{Theorem \ref{teo:average}} we get
$$
{\kappa_{av}}((A,b),y)=\frac{\|A^{-1}\|_F\cdot\sqrt{1+\|y\|^2}}{\sqrt{n^2+n}}
\leq\frac{\kappa((A,b),y)}{\sqrt{n+1}}.
$$

For the $k$th-componentwise condition number, we have that
$$
{\kappa_{av}}((A,b),y,k)=\left[\frac{1}{\sqrt{\pi}}\cdot
\frac{\Gamma\left(\frac{n^2}{2} \right)}{\Gamma\left(\frac{n^2+p}{2} \right)}\cdot\Gamma\left(\frac{p+1}{2}\right)\right]^{1/p}\cdot
\kappa((A,b),y).
$$
A similar result was proved in \cite{S-W-W-W}.

In \cite{Edelman}, it is proved that the expected value of the relative condition number $\kappa_{rel}(A)=\|A\|\cdot\|A^{-1}\|$ of a random matrix $A$ whose elements are i.i.d standard normal, satisfies:
$$
\mathbb{E}(\log \kappa_{rel}(A))= \log n +c +o(1),
$$
as $m\to \infty$, where $c\approx 1.537$.
If we consider the relative average condition number defined in (\ref{def:cnrav}),  we get from (\ref{eq:acnlc})
$$
\mathbb{E}(\log {\kappa_{rel}}_{av}(A))= \frac12\log n +c+o(1),
$$
as $m\to \infty$.


\vspace{10pt}

\subsection{Eigenvalue and Eigenvector Problem}
Let $X= \mnc$ be the space of $n\times n$ complex matrices with the Frobenius Hermitian inner product,
$Y=\mathbb{P}(\C^{n+1})\times\C$ and $V=\{(A,v,\lambda):\, Av=\lambda v \}$.
Then for $(A,v,\lambda)\in V\setminus\Sigma$  the condition matrices $DG_1$ and $DG_2$ associated with the eigenvector and eigenvalue problem are
$$
DG_1(A)\dot A=\left(\pi_{v^\perp}(\lambda I-A)|_{v^\perp}  \right)^{-1}\left(\pi_{v^\perp }\dot A v\right)\quad \mbox{and}
\quad DG_2(A)\dot A= \frac{\pes{u}{\dot A v}}{\pes{u}{v}},
$$
where $u$ is some left eigenvector associated  with $\lambda$, i.e. $u^*A=\overline{\lambda}u^*$ (see [Bez IV]).
The associated condition numbers are:
\beq\label{eq:condvalp}
\kappa_1(A,v)= \left\| \left(\pi_{v^\perp}(\lambda I-A)|_{v^\perp}  \right)^{-1} \right\|
\quad\mbox{and}\quad
\kappa_2(A,\lambda)=\frac{\|u\|\cdot\|v\|}{|\pes uv|}.
\eeq
From our \textit{Theorem \ref{teo:average}}, we get the respective average condition numbers:
$$
{\kappa_{av}}_1(A,v)= \frac 1n \left\| \left(\pi_{v^\perp}(\lambda I-A)|_{v^\perp}  \right)^{-1} \right\|_F\leq \frac{1}{\sqrt n}\kappa_1(A,v),
$$
$$
{\kappa_{av}}_2(A,\lambda)=\frac1n \kappa_2(A,\lambda).
$$
A similar result for ${\kappa_{av}}_2(A,\lambda)$ was proved in \cite{Stewart}.

\vspace{10pt}


\subsection{Finding Kernels of Linear Transformations}

Let $\mmnc$ be the linear space of $k\times p$ complex matrices with the Frobenius Hermitian inner product, i.e. $\pes{A}{B}_F=\mbox{trace}(B^*A)$ ($B^*$ is the adjoint of $B$),  and $\mathcal{R}_r \subset\mmnc$ the subset of matrices of rank $r$.
We consider the problem of solving the system of linear equations $Ax=0$.
For this purpose, we introduce the  \textit{Grassmannian}  $\mathbb{G}_{p,\ell}$ of complex subspaces of dimension $\ell$ in $\C^p$, where $\ell=\dim (\ker A)=p-r$ .

The input space $X=\mathcal{R}_r$ is a smooth submanifold of $\mmnc$ of complex dimension $(k+p)r-r^2$ (see \cite{Dedieu}). Thus, it has a natural Hermitian structure induced by the Frobenius metric on $\mmnc$.

 In what follows, we identify $\mathbb{G}_{p,\ell}$ with the quotient $\mathbb{S}_{p,\ell}/\mathcal{U}_\ell$ of the \textit{Stiefel} manifold
 $$\mathbb{S}_{p,\ell}:=\{M\in\mathscr{M}_{p,\ell}(\C):\,M^*M=I \}$$
  by the unitary group $\mathcal{U}_\ell\subset \mathscr{M}_\ell(\C)$, which acts on the right of  $\mathbb{S}_{p,\ell}$ in the natural way (see \cite{Dedieu}).
  Then, the complex dimension of the output space $Y=\mathbb{G}_{p,\ell}$ is $(p-r)r$.

We will use the same symbol to represent an element of $\mathbb{S}_{p,\ell}$ and its class in $\mathbb{G}_{p,\ell}$.
 The manifold $\mathbb{S}_{p,\ell}$  has a canonical Hermitian structure induced by the Frobenius norm in $\mathscr{M}_\ell(\C)$.
On the other hand, $\mathcal{U}_\ell$ is a Lie group of isometries acting on $\mathbb{S}_{p,\ell}$. Therefore,  $\mathbb{G}_{p,\ell}$ is a Homogeneous space (see \cite{G-H-L}), with a natural Riemannian structure that makes the projection $\mathbb{S}_{p,\ell}\rightarrow\mathbb{G}_{p,\ell}$ a Riemannian submersion.
The orbit of $M\in\mathbb{S}_{p,\ell}$ under the action of the unitary group $\mathcal{U}_\ell$, namely, $\mbox{o}_\ell(M)=\{MU:\,U\in \mathcal{U}_\ell\}$, defines a smooth submanifold of $\mathbb{S}_{p,\ell}$.  In this form we can define a local chart of a small neighborhood of $M\in\mathbb{G}_{p,\ell}$ from the affine spaces
$$
M+ T_M\mbox{o}_\ell(M)^{\perp}
$$
where $T_M\mbox{o}_\ell(M)^{\perp}$ is the orthogonal complement of $T_M\mbox{o}_\ell(M)$  in $T_M \mathbb{S}_{p,\ell}$.
Implicit differentiation in local coordinates of the input-output map $G$ at the point $G(A)=M\in\mathbb{S}_{p,\ell}$  yields
\beq\label{eq:cmkp}
DG(A)(\dot A)=-A^\dagger \dot A M,\
\eeq
where $\dot A\in T_A \mathcal{R}_r $, $A^\dagger$ is the Moore-Penrose inverse of $A$ and $A^\dagger \dot A M\in T_M\mbox{o}_\ell(M)^{\perp}$.

One way to compute the singular values of the condition matrix described in (\ref{eq:cmkp}), is to take an orthonormal basis in $\mmnc$  which diagonalizes $A$.
 From the singular value decomposition, there exists orthonormal basis $\{u_1,\ldots,u_k\}$ of $\C^k$, and $\{v_1,\ldots,v_p\}$ of $\C^p$,  such that $Av_i=\sigma_i u_i$ for $i=1,\ldots,r$, and $Av_i=0$ for $i=r+1,\ldots,p$.
 Thus, $\{u_i v_j^*:\,i=1,\ldots,k;\,j=1,\ldots,p \}$ is an orthonormal basis  of $\mmnc$ which diagonalizes $A$.
 In this basis the tangent space $T_A \mathcal{R}_r$ is the orthogonal complement of the subspace generated by $\{u_i v_j^*:\,i=r+1,\ldots,k;\,j=r+1,\ldots,p \}$.
From where we conclude that
$$
{\kappa}(A,M)=\|DG(A)\|=\|A^\dagger\|,\quad
{\kappa}_F(A,M)=\sqrt{p-r}\cdot \|A^\dagger\|_F.
$$
From our \textit{Theorem \ref{teo:average}},
$$
{\kappa_{av}}(A,M)=\frac{\sqrt{p-r}}{\sqrt{(k+p-r)r}}\cdot\|A^\dagger\|_F\leq  \sqrt{\frac{p(p-r)}{(k+p-r)r}}\cdot   \kappa(A,M).
$$

In \cite{Beltran}, it is proved that
$$
\mathbb{E}(\log \kappa_{rel}(A):\, A\in \mathcal{R}_r)\leq\log\left[\frac{k+p-r}{k+p-2r+1}\right] +2.6,
$$
where the expected value is computed with respect to the normalized naturally induced measure in $\mathcal{R}_r$. Our \textit{Theorem \ref{teo:average}} immediately yields a bound for the average relative condition number, namely,
$$
\mathbb{E}(\log {\kappa_{rel}}_{{av}}(A):\, A\in \mathcal{R}_r) \leq \frac{1}{2}\log \left[\frac{(k+p-r)r}{(k+p-2r+1)^2p(p-r)}\right] +2.6.
$$

\vspace{10pt}


\subsection{Finding Roots Problem I: Univariate Polynomials}
We start with the case of one polynomial in one complex variable.
Let $X=\mathcal{P}_d=\{f:\,f(z)=\sum_{i=0}^d f_iz^i, \,f_i\in\C\}$. Identifying $\mathcal P_d$
with $\C^{d+1}$ we can define two standard inner products in the space $\mathcal P_d$:\\
- Weyl inner product:
\beq\label{def:Weyl}
\pes{f}{g}_W:=\sum_{i=0}^d f_i\overline{g_i}\binom{d}{i}^{-1};
\eeq
- Canonical Hermitian inner product:
\beq\label{def:canonical}
\pes{f}{g}_{\C^{d+1}}:=\sum_{i=0}^d f_i\overline{g_i}.
\eeq
The solution variety is given by $V=\{(f,z):\,f(z)=0 \}$. Thus, by implicit differentiation
$$
DG(f)(\dot f)=-\left(f'(\zeta)\right)^{-1}\dot{f}(\zeta).
$$
We denote by $\kappa_{W}$ and $\kappa_{\C^{d+1}}$  the condition numbers with respect to the Weyl
and Euclidean inner product.
The reader may check that
$$
\kappa_{W}(f,\zeta)=\frac{(1+|\zeta|^2)^{d/2}}{|f'(\zeta)|}	\quad\mbox{and}\quad
\kappa_{\C^{d+1}}(f,\zeta)=\frac{\sqrt{\sum_{i=0}^d|\zeta|^{2i}}}{|f'(\zeta)|},
$$
(for a proof see \cite{B-C-S-S}, p. 228 ).
From \textit{Theorem \ref{teo:average}}, we get:
$$
{\kappa_{av}}_{W}^{[2]}(f,\zeta)=\frac{1}{\sqrt {2(d+1)}}\kappa_{W}(f,\zeta)
,	\qquad
{\kappa_{av}}_{\C^{d+1}}^{[2]}(f,\zeta)=\frac{1}{\sqrt{2(d+1)}}\kappa_{\C^{d+1}}(f,\zeta).
$$

\vspace{8pt}

\subsection{Finding Roots Problem II: Systems of Polynomial Equations}

We now study the case of complex homogeneous polynomial systems.
Let $\mathcal{H}_{(d)}$ the space of systems $f:\C^{n+1}\to\C^n$, $f=(f_1,\ldots,f_n)$ where each $f_i$ is a homogenous polynomial
 of degree $d_i$. We consider $\mathcal{H}_{(d)}$ with the homogeneous analogous of the
 Weyl structure defined above (see Chapter 12 of \cite{B-C-S-S}  for details).

\vspace{5pt}

Let $X=\mathbb{P}(\mathcal{H}_{(d)})$ and $Y=\mathbb{P}(\C^{n+1})$
and $V=\{(f,\zeta):\,f(\zeta)=0 \}$. We denote by $N=\sum_{i=1}^n \binom{d_i+n}{n}-1$
 the complex dimension of $X$. We may think of $2N$ as the size of the input.

Then,
$$
DG(f)(\dot f)=-\left(Df(\zeta){_{\zeta^\perp}}\right)^{-1}\dot f(\zeta),
$$
and the condition number is
$$
\kappa_W(f,\zeta)= \left\| \left(Df(\zeta){_{\zeta^\perp}}\right)^{-1}
\right\|,
$$
where some norm $1$ affine representatives of $f$ and $\zeta$ have been chosen (cf.\cite{B-C-S-S}). Associated
with this quantity, we consider
\beq\label{def:cnf}
\kappa_W(f):=\sqrt{\frac{1}{\mathcal{D}}
\sum_{\{\zeta:\,f(\zeta)=0 \}}\kappa_W(f,\zeta)^2 },
\eeq
where $\mathcal D=d_1\cdots d_n$ is the number of projective solutions of a generic system.

The expected value of $\kappa_W^2(f)$ is an essential ingredient in the complexity analysis
of path-following methods (cf. \cite{BezIV}, \cite{B-P}).
In \cite{B-P} the authors proved that
\beq\label{eq:BP}
\mathbb{E}_{f}\left[\kappa_W(f))^2\right]
\leq 8 n N,
\eeq
where $f$ is chosen at random with the Weyl distribution.

The relation between complexity theory and ${\kappa_{av}}$ is not clear yet.
However, it is interesting to study the expected value of the ${\kappa_{av}}$-analogous of
equation (\ref{eq:BP}), namely
$$
{\kappa_{av}}_W(f):=
\sqrt{\frac{1}{\mathcal{D}}
\sum_{\{\zeta:\,f(\zeta)=0 \}}{\kappa_{av}}_W(f,\zeta)^2 }.
$$
From our \textit{Theorem \ref{teo:average}} we get,
$$
{\kappa_{av}}_W(f,\zeta)\leq\frac{\kappa_W(f,\zeta)}{\sqrt{N/n}},
\qquad
\mathbb{E}_{f}\left[{\kappa_{av}}_W(f))^2\right]\leq 8\,n^2.
$$
Note that the last bound depends on the number of unknowns $n$, and not on the size of the
 input $n\ll N$.


\vspace{10pt}

\section{Proof of the main Theorem }

In the case of complex manifolds, the condition matrix turns to be an $n\times n$ complex matrix.
In what follows, we identify it with the associated $2n\times2n$ real matrix. We center our attention in the real case.

The main theorem  follows immediately  from \textit{Lemma \ref{lem:unigaus}} and
\textit{Proposition \ref{prop:avergaus}} below.

\begin{lem}\label{lem:unigaus}
Let $\eta$ be a Gaussian standard random vector in $\R^m$. Then
$$
{\kappa_{av}}^{[p]}(x,y)=
\frac{1}{\sqrt 2}\left[\frac{\Gamma\left(\frac{m}{2} \right)}{\Gamma\left(\frac{m+p}{2} \right)}\right]^{\frac 1p}
\cdot\left[\mathbb{E}(\|DG(x)\eta\|^p)  \right]^{1/p},
$$
where $\mathbb{E}$ is the expectation operator and $\|\cdot\|$ is the Euclidean norm in $\R^n$.

\end{lem}
\begin{proof}
Let $f:\R^m\to\R$ be the continuous function given by
$$
f(v)=\|DG(x)v\|.
$$
Then
$$
\left[\mathbb{E}(\|DG(x)\eta\|^p)  \right]^{1/p}=\left[\frac{1}{(2\pi)^{m/2}}\int_{\R^m} f(v)^p\cdot e^{-\|v\|^2/2}\,dv  \right]^{1/p}.
$$
Integrating in polar coordinates we get that
\beq\label{eq:intpol}
\mathbb{E}(\|DG(x)\eta\|^p) =\frac{I_{m+p-1}}{(2\pi)^{m/2}}\cdot   \int_{S^{m-1}} f^p\,dS^{m-1}
\eeq
where
$$
I_{j}=\int_0^{+\infty}\rho^{j}\,e^{-\rho^2/2}\,d\rho,\quad j\in\N.
$$
Making the change of variable $u=\rho^2/2$ we obtain
$$
I_j=
2^{\frac{j-1}{2}}\Gamma(\frac{j+1}{2}),
$$
therefore
\beq\label{eq:ij}
I_{m+p-1}
=2^{\frac{m+p-2}{2}}\cdot\Gamma\left(\frac{m+p}{2}\right).
\eeq
Then joining together (\ref{eq:intpol}) and (\ref{eq:ij}) we obtain
the result.
\end{proof}

\begin{prop}\label{prop:avergaus}
$$
\mathbb{E}(\|DG(x)\eta\|^p)  =\mathbb{E}(\|\eta_{\sigma_1,\ldots,\sigma_n}\|^p)
$$
where
$\eta_{\sigma_1,\ldots,\sigma_n}$ is a centered Gaussian vector in $\R^n$ with diagonal covariance matrix
$Diag(\sigma_1^2,\ldots,\sigma_n^2)$.
\end{prop}
\begin{proof}
Let $DG(x)=UDV$ be a singular value decomposition of $DG(x)$.
By the invariance of the Gaussian distribution under the action of the orthogonal group in $X$,
$V\eta$ is again a standard Gaussian random vector. Then,
$$
\mathbb{E}(\|DG(x)\eta\|^p) =
\mathbb{E}(\|UD\eta\|^p),
$$
and by the invariance under the action of the orthogonal group of the Euclidean norm, we get
$$
\mathbb{E}(\|DG(x)\eta\|^p) =
\mathbb{E}(\|D\eta\|^p).
$$
Finally $D\eta$ is a centered Gaussian vector in $\R^n$ with covariance matrix
$Diag(\sigma_1^2,\ldots,\sigma_n^2)$, and the proposition follows.
For the case $p=2$
$$
{\kappa_{av}}(x,y)=
\left[\mathbb{E}\left(\sigma_1^2\eta_1^2+\ldots+\sigma_n^2\eta_n^2\right)\right]^{1/2}
$$
where $\eta_1,\ldots,\eta_n$ are i.i.d. standard normal in $\R$.
Then
$$
{\kappa_{av}}(x,y)=\left( \sum_{i=1}^n\sigma_i^2\right)^{1/2}=\kappa_F(x,y).
$$
\end{proof}

\vspace{10pt}

\vspace{8pt}

\textbf{Acknowledgments}
\\The author thanks Professors Carlos Beltr\'an, Michael Shub and Mario Wschebor for useful discussions.

\vspace{10pt}



\begin{thebibliography}{99}
\bibitem{Beltran}C. Beltr\'an: Estimates On The Condition Number Of Random,
Rank-Deficient Matrices. (to appear)
\bibitem{B-P} C. Beltr\'an and L.M. Pardo: Fast linear homotopy to find approximate zeros of polynomial systems. (to appear)
\bibitem{B-C-S-S} L. Blum, F. Cucker, M. Shub and S. Smale:
Complexity and real computation.
\emph{ Springer-Verlag}, New York, 1998.
\bibitem{BCL} P. B\"urgisser, F. Cucker and M. Lotz: Smoothed analysis of complex
conic condition numbers, \emph{J. Math. Pures et Appl.} 86, 293-309, 2006.
\bibitem{Bur} P. B\"urgisser:
Smoothed Analysis of Condition Numbers, (to appear).
\bibitem{Dedieu} J.P. Dedieu:
Points Fixes, Z\'eros et la M\'ethode de Newton. Math\'ematiques \& Applications \emph{SMAI-Springer Verlag },  54, 2006.
\bibitem{Edelman} A. Edelman: Eigenvalues and condition numbers of random matrices, \emph{SIAM J. of Matrix Anal.
and Applic.}, 9, 543-560, 1988.
\bibitem{G-H-L} S. Gallot, D. Hulin and J. Lafontaine:
Riemannian geometry. Third edition. Universitext. \emph{Springer-Verlag}, Berlin, 2004.
\bibitem{BezIV} M. Shub and S. Smale:
Complexity of Bezout's theorem. IV. Probability of success; extensions.
\emph{SIAM J. Numer. Anal.}  33 , no. 1, 128-148, 1996.
\bibitem{Smale} S. Smale: On the eficiency of Algorithms of analysis, \emph{Bull. Amer. Math. Soc.}, 13 ,
 87-121, 1985.
\bibitem{SpTeng} D.A. Spielman and S.H. Teng: Smoothed analysis of algorithms, \emph{ICM 2002}, Beijing, Vol. I, pp. 597-606, 2002.
\bibitem{Stewart} G. W. Stewart: Stochastic perturbation theory,
\emph{ SIAM Rev.}  32  ,  no. 4, 579-610, 1990.
\bibitem{Tref} L.N. Trefethen and D. Bau, III: Numerical linear algebra, \emph{Society
for Industrial and Applied Mathematics (SIAM)}, Philadelphia, PA, 1997.
\bibitem{S-W-W-W} N. Weiss, G.W. Wasilkowski, H. Wozniakowski,  and  M. Shub:
Average condition number for solving linear equations.
\emph{Linear Algebra Appl.} 83,  79-102, 1986.
\bibitem{Wschebor} M. Wschebor: Smoothed Analysis of $\kappa(a)$, \emph{J. of Complexity}, 20, 97-107, 2004.
\end{thebibliography}
\end{document}